\documentclass[12pt]{article}
\usepackage{amssymb}
\usepackage{graphicx}
\usepackage{latexsym}
\def\part#1{\frac{\partial\phantom{q}}{\partial#1}}

\newenvironment{rmk}{\begin{trivlist}\item[]{\bf Remark:} }
{\end{trivlist}}
\newenvironment{rmks}{\begin{trivlist}\item[]{\bf Remarks:} }
{\end{trivlist}}
\newenvironment{ex}{\begin{trivlist}\item[]{\bf Example:} }
{\end{trivlist}}
\newenvironment{prf}{\begin{trivlist}\item[]{\bf Proof:} }
{\hfill $\Box$ \end{trivlist}}

\newtheorem{thm}{Theorem}

\newtheorem{prp}[thm]{Proposition}

\def\End{\mathop{\rm End}\nolimits}
\def\Hom{\mathop{\rm Hom}\nolimits}

\def\adj{\mathop{\rm adj}\nolimits}

\newcommand{\R}{\mathbf{R}}
\newcommand{\C}{\mathbf{C}}

\newcommand{\PP}{{\mathbf {\rm P}}}

\begin{document}
\title{Poisson modules and generalized geometry}
 \author{Nigel Hitchin\\[5pt]}
\maketitle
\centerline{\it Dedicated to Shing-Tung Yau on the occasion of his 60th birthday}

\section{Introduction}

Generalized complex structures were introduced as a common format for discussing both symplectic and complex manifolds, but the most interesting examples are hybrid objects -- part symplectic and part complex. One such class of examples  consists of holomorphic Poisson surfaces, but in \cite{CG1},\cite{CG2} Cavalcanti and Gualtieri also construct generalized complex 4-manifolds with similar features which are globally neither complex nor symplectic.

In Gualtieri's development of the subject \cite{G0},\cite{G1} he introduced generalized analogues of a number of  familiar concepts in complex geometry, and notably the idea of a generalized holomorphic bundle. In the symplectic case this is simply a flat connection, but in the Poisson case it is more interesting and coincides with the notion of Poisson module: a locally free sheaf with a Poisson action of the sheaf of functions on it. These play a significant role in Poisson geometry, and can be thought of as semiclassical limits of noncommutative bimodules.

A Poisson structure on a complex surface is a section $\sigma$ of the anticanonical bundle. It vanishes in general on an elliptic curve. We begin  this paper by using algebraic geometric methods to  construct rank two
Poisson modules on such a surface. The data for this particular construction is located on the elliptic curve: a line bundle together with a pair of sections with no common zero.

Now if we consider the surface as a generalized complex manifold, then where $\sigma\ne 0$, $\sigma^{-1}=B+i\omega$ and we regard the generalized complex structure as being defined by the symplectic  form $\omega$ transformed by a B-field $B$.  Where $\sigma=0$  it is the transform of a complex  structure.  The nonholomorphic examples in \cite{CG1},\cite{CG2}  also have a 2-torus on which the generalized complex structure changes type from symplectic to complex. Moreover the torus acquires the structure of an elliptic curve. This provokes the natural question of whether a holomorphic line bundle on this curve with a  pair of sections  can generate a generalized holomorphic bundle in analogy with the holomorphic Poisson case. This we answer in the rest of the paper, and in fact conclude from the general result that a line bundle with a single section is sufficient.  The proof entails replacing the algebraic geometry of the Serre construction by a differential geometric version, and using this as a model in the generalized case.

As an application, we adapt a construction of Polishchuk in \cite{P} to define generalized complex structures on $\PP^1$-bundles over the examples of Cavalcanti and Gualtieri.

\section{Poisson modules}
\subsection{Definitions}
Let $M$ be a holomorphic Poisson manifold, defined by a section $\sigma$ of $\Lambda^2T$, then the Poisson bracket of two locally defined holomorphic functions $f,g$ is $\{f,g\}=\sigma(df,dg)$. Algebraically, a Poisson module is a locally free sheaf ${\mathcal O}(V)$ with an action $s\mapsto \{f,s\}$ of the structure sheaf with the properties
\begin{itemize}
\item
$\{f,gs\}=\{f,g\}s+g\{f,s\}$
\item
$\{\{f,g\},s\}=\{f,\{g,s\}\}-\{g,\{f,s\}\}.$
\end{itemize}
The first equation defines a first order linear differential operator 
$$D:{\mathcal O}(V)\rightarrow {\mathcal O}(V\otimes T)$$
where $\{f,s\}=\langle Ds,df \rangle$. This is simply a holomorphic differential operator whose symbol is $1\otimes\sigma: V\otimes T^*\rightarrow V\otimes T$. 

The second equation is a zero curvature condition. Relative to a local basis $s_i$ of $V$, $D$ is defined by a ``connection matrix" $A$ of vector fields:
$$Ds_i=\sum_j s_j\otimes A_{ji}$$
and the  condition  becomes ${\mathcal L}_A\sigma=A^2\in \End(V\otimes \Lambda^2T)$. When $\sigma$ is non-degenerate and identifies $T$ with $T^*$ then this is a flat connection.

\begin{ex} If $X=\sigma(df)$ is the Hamiltonian vector field of $f$ then the Lie derivative ${\mathcal L}_X$ acts on tensors  but the action in general  involves the second derivative of $f$. However for the canonical line bundle $K=\Lambda^nT^*$ we have 
$${\mathcal L}_X(dz_1\wedge\dots\wedge dz_n)=\sum_i \frac{\partial X_i}{\partial z_i} (dz_1\wedge\dots\wedge dz_n)$$
and
$$\sum_i\frac{\partial X_i}{\partial z_i}=\sum_{i,j}\frac{\partial}{\partial z_i}(\sigma^{ij}\frac{\partial f}{\partial z_j})=\sum_{i,j}\frac{\partial \sigma^{ij}}{\partial z_i}\frac{\partial f}{\partial z_j}$$
which involves only the first  derivative of $f$. Thus 
$$\{f,s\}={\mathcal L}_Xs=\langle Ds, df\rangle$$
defines a first order operator. The second condition follows from the integrability of the Poisson structure: since $\sigma(df)=X, \sigma(dg)=Y$ implies $\sigma(d\{f,g\})=[X,Y]$, it follows that 
$$\{\{f,g\},s\}={\mathcal L}_{[X,Y]}s=[{\mathcal L}_X,{\mathcal L}_Y]s=\{f,\{g,s\}\}-\{g,\{f,s\}\}.$$
This clearly holds for any power $K^m$. 
\end{ex}
\subsection{A construction}\label{con1}
If  a rank $m$ vector bundle $V$ is a Poisson module, then so is the line bundle $\Lambda^mV$.
Now let $M$ be a complex surface and $V$ a rank 2 holomorphic vector bundle with $\Lambda^2V\cong K^*$ (a line bundle which, as noted above, is a Poisson module for any Poisson structure). Suppose $V$ has two sections $s_1,s_2$ which are generically linearly independent. Then $s_1\wedge s_2$ is a holomorphic section $\sigma$ of $\Lambda^2V\cong K^*$ and so defines a Poisson structure. It vanishes on a curve $C$. Moreover, where $\sigma\ne 0$, $s_1,s_2$ are linearly independent  and define a trivialization of $V$ and hence a flat connection.
\begin{prp} The flat connection extends to a Poisson module structure on $V$.
\end{prp}
\begin{prf}  Let $(u_1,u_2)$ be a local holomorphic basis for $V$ in a neighbourhood of a point of $C$, then 
$$s_i=\sum_j P_{ji} u_j$$
and $\det P=0$ is the local equation for $C$.  Now a connection matrix for $D$ with $Ds_i=0$ in the basis $(u_1,u_2)$ is given by a matrix $A$ of vector fields such that
$$0=Ds_i=D(\sum_j P_{ji} u_j)=\sum_j\sigma(dP_{ji})u_j+\sum_{jk}P_{ji}A_{kj}u_k$$
which has solution
$$A=-\sigma(dP)P^{-1}=-\sigma(dP)\frac{\adj P}{\det P}.$$
This is smooth since $\det P$  divides  $\sigma$.
\end{prf}

Finding rank 2 bundles with two sections is a problem that has been considered before, most notably in the study of charge $2$ instanton bundles on $\PP^3$ \cite{H}. Our case is similar but one dimension lower.
The choice of two sections defines an extension of sheaves:
$$0\rightarrow {\mathcal O}\oplus {\mathcal O}\rightarrow {\mathcal O}(V)\rightarrow {\mathcal O}_C(L^*K_M^*)\rightarrow 0$$
where $L$ is the line bundle on the elliptic curve $C$ where the two, now linearly dependent, sections take their value.
Such an extension  is classified by an element of 
global $ {\mathrm{Ext}}^1({\mathcal O}_C(L^*K_M^*),{\mathcal O})\otimes \C^2$
but local duality gives an isomorphism of sheaves  
${\mathcal Ext}^1({\mathcal O}_C,K_M)\cong {\mathcal O}_C(K_C).$ 
Hence
$${\mathrm{Ext}}^1({\mathcal O}_C(L^*K_M^*),{\mathcal O})\otimes \C^2\cong H^0(C,L)\otimes \C^2$$
and we are looking for a pair   of sections of the line bundle $L$ on $C$. To get a locally free sheaf  we need the pair to have no common zeros. 

To relate this data to the vector bundle, note that any one of the sections is non-vanishing outside $C$ and so the number of zeros is the same as the number of zeros of a section of $L$. Counting multiplicities, this means that 
$$c_1(L)=c_2(V).$$

Another way of recording the information is to consider the meromorphic function $s_1/s_2$ on $C$.

\section{The Serre construction}
\subsection{The algebraic approach}
To obtain an analogue of the  construction above when $M$ is generalized complex, we have to replace 
the sheaf theory by a more analytic method. To prepare for this we consider next  the Serre construction of rank 2 holomorphic vector bundles. This is the question of constructing a vector bundle with at least {\it one} section and not two as  above. For surfaces, a good source is Griffiths and Harris \cite{GH}. We begin with a sheaf-theoretic approach.

The problem is this: given $k$ points $x_i\in M$, find a rank 2 vector bundle $V$ with a section $s$ which vanishes non-degenerately at the $k$ points. We shall assume, as previously, that $\Lambda^2V\cong K^*$. The derivative of $s$ at a zero $x$ is a well defined element of $(V\otimes T^*)_x$ and we can take
$$\det ds(x)\in (\Lambda^2V\otimes \Lambda^2T^*)_x$$
which is canonically the complex numbers  thanks to the isomorphism $\Lambda^2V\cong K^*$. It is non-zero since the zero is nondegenerate. So at each zero $x_i$, $s$ has a residue $(\det ds(x_i))^{-1}=\lambda_i\ne 0$. Then, as in \cite{GH} Chapter 5:

\begin{prp} \label{Serre} Given k points $X= \{x_1, . . . , x_k\} \in M$,
and $\lambda_i\in \C^*$  such that $\lambda_1+\dots +\lambda_k = 0$, there exists a rank 2 vector
bundle V with $\Lambda^2V = K^*$  and a section s such that $s(x_i) = 0$ 
and $\det ds(x_i) =\lambda^{-1}_i$.
\end{prp}

In this case the bundle is given by  an extension of sheaves:
$$0\rightarrow {\mathcal O}\rightarrow {\mathcal O}(V)\rightarrow {\mathcal I}_X\otimes K^*\rightarrow 0.$$
To link things up with the previous section we ask when there is a second section. The exact cohomology sequence gives:
$$0\rightarrow H^0(M,{\mathcal O})\rightarrow H^0(M,{\mathcal O}(V))\rightarrow H^0(M,{\mathcal I}_X\otimes K^*)\rightarrow H^1(M,{\mathcal O})\rightarrow$$
Now if the points lie on the zero set of a Poisson structure $\sigma$, then $\sigma$ lies in the space $H^0(M,{\mathcal I}_X\otimes K^*)$.  If $H^1(M,{\mathcal O})=0$ then $\sigma$ pulls back to a second section.

\begin{rmks}
\noindent 1. Our assumption that $\sigma$ vanishes nondegenerately on a single elliptic curve actually implies that  $H^1(M,{\mathcal O})=0$ when the surface is algebraic. This follows from the classification 
of \cite{BM}.

\noindent 2. In the previous section, the data for the construction of the bundle was a pair of sections $s_1,s_2$ on $C$. One might ask where the $\lambda_i$ come from if we choose just one of these, $s_1$. In fact $d\sigma$ restricted to $C$ gives an isomorphism between  the normal bundle $N$ and $K^*$. But $K^*\cong NK^*_C$ so we get a canonical non-vanishing vector field on $C$ (the so-called {\it modular vector field}). Its inverse is a nonvanishing differential $\alpha$ and then we can take $\lambda_i$ to be the residue of the differential $(s_2/s_1)\alpha$ at $x_i\in C$. The sum of the residues of a differential on a curve is of course zero.
\end{rmks}
\subsection{The analytical approach}
We now reformulate the Serre construction in Dolbeault terms (see also Chapter 10.2 in \cite{DK}). First consider the sequence of sheaves:
$$0\rightarrow {\mathcal O}\rightarrow {\mathcal O}(V)\stackrel{\pi}\rightarrow {\mathcal I}_X\otimes K^*\rightarrow 0.$$
This is an extension of line bundles outside the points $x_i$, and the standard way  to obtain a Dolbeault representative for this is to choose a Hermitian metric on $V$, and form  the orthogonal complement of the trivial subbundle.  Restrict  $\bar\partial$ to this line bundle.  Since the homomorphism $\pi$ in the complex is holomorphic, we obtain a $\bar\partial$-closed $(0,1)$-form with values in $\Hom(K^*,{\mathcal O})\cong K$. In other words a $(2,1)$-form.

In our case this will acquire singularities at the points $x_i$, so let us consider the local model.
Let   $x_i$ be given by the origin in coordinates $z_1,z_2$, then the two maps in the complex are represented by 
$1\mapsto (z_1,z_2)$ and $(u,v)\mapsto -z_2u+z_1v.$ Using the trivial Hermitian structure we take the orthogonal complement of $(z_1,z_2)$ to give 
 $$\bar\partial\left(\frac{1}{r^2}(-\bar z_2,\bar z_1)\right)=\frac{1}{r^4}(\bar z_2d\bar z_1-\bar z_1d\bar z_2)(z_1,z_2)$$
and then
\begin{equation}
A_0=\frac{1}{r^4}dz_1\wedge dz_2\wedge (\bar z_2d\bar z_1-\bar z_1d\bar z_2)
\label{fund}
\end{equation}
is the required $(2,1)$ form.

Now, using the flat metric on $\C^2$, we calculate 
$$\ast d\left (\frac{1}{r^2}\right)^{2,1}=\frac{1}{4r^4}dz_1\wedge dz_2\wedge (\bar z_2d\bar z_1-\bar z_1d\bar z_2)$$
and in four dimensions $1/r^2$ is, up to a universal constant, the fundamental solution of the Laplacian. It follows that 
$$\bar\partial \left(\frac{1}{r^4}dz_1\wedge dz_2\wedge (\bar z_2d\bar z_1-\bar z_1d\bar z_2)\right)= c \delta(0).$$

Thus, a distributional $(2,1)$ form $A$ which has the form (\ref{fund}) at each point $x_i$ and is smooth elsewhere, is an analytical way of defining the vector bundle $V$ with a section.  

\begin{rmk} Near $x_i$ we have a non-holomorphic basis for $V$ defined by $(z_1,z_2),(-\bar z_2,\bar z_1)$. This is obtained from a  trivialization which extends to $x_i$ by the gauge transformation on $\C^2\setminus \{0\}$ 
$$\pmatrix  {z_1 & z_2\cr
                           -\bar z_2 & \bar z_1}.$$
For example the flat trivialization of the spinor bundle on $\R^4$, the complement of a point in $S^4$, extends this way.
\end{rmk}

The proof of Proposition  \ref{Serre} now goes as follows. Consider the distributional form (or {\it current}) $T$ defined by taking delta functions at the points $x_i$: 
$$T=\sum_i\lambda_i\delta(x_i)\in \Omega^{2,2}.$$
This defines a class in $H^2(M,K)$. Since $H^2(M,K)$ is dual to $H^0(M,{\mathcal O})\cong \C$, this class is determined by evaluating it on the function $1$. But
$$\langle T, 1\rangle =\sum_i\lambda_i$$
so if the $\lambda_i$ sum to zero the class is zero. 

Now choose a Hermitian metric on $M$, flat near the $x_i$. From harmonic theory  there is a current $S\in \Omega^{2,2}$ such that $\bar\partial\bar\partial^*S=T$ and then we take
$$A=\bar\partial^*S$$
to be the distribution defining the extension. Near $x_i$,
$$\bar\partial\bar\partial^*(S-k\ast \!\frac{1}{r^2})=0$$
so by elliptic regularity the difference is smooth and $A$ defines a holomorphic structure on the bundle obtained as in the Remark above.

\begin{rmk} The 't Hooft construction of $SU(2)$ instantons on $\R^4=\C^2$ defines an anti-self-dual connection (and a fortiori a holomorphic structure) from a harmonic function of the form
$$\phi=\sum_i\frac{1}{\vert {\mathbf x}- {\mathbf x_i}\vert^2}.$$
Its twistor interpretation is the Serre construction for the lines in $\PP^3$ corresponding to the points $x_i\in \R^4$ (see \cite{AW}). This is a model for the above reformulation.

\end{rmk}

\subsection{The second section}
We now look analytically for a second section of the vector bundle $V$. The bundle outside of $X$ is  a direct sum $1\oplus K^*$ with $\bar\partial$-operator defined by 
$$\bar\partial (u,v)=(\bar\partial u+Av,\bar\partial v).$$
In this description, the section coming from the Serre construction is $s_1=(1,0)$, and we want a second one $s_2$ such that $s_1\wedge s_2=\sigma$, so we write  $s_2=(u,\sigma)$. For holomorphicity we need 
$$\bar\partial u+ A\sigma=0.$$
Now
 $$T=\sum_i\lambda_i \delta(x_i)\in \Omega^{2,2}$$
so consider 
$\sigma T\in \Omega^{0,2}$. Evaluating $\sigma T$ on a $(2,0)$ form $\alpha$ gives
$$\sum_i\lambda_i \sigma(x_i)\alpha(x_i)=0$$
if the points $x_i$ lie on the zero set of $\sigma$. Thus $\sigma T =0$.

We have $\bar\partial A=T$ and $\sigma$ is holomorphic so $\bar \partial (A\sigma)=0$. But if 
$H^1(M,{\mathcal O})=0$ then this implies that $A\sigma=-\bar\partial u$ for the required distributional section $u$.

\section{Generalized geometry}
\subsection{Basic features}
Before we  adapt this method to generalized complex manifolds, we  review here the basic features. For more details see \cite{G0},\cite{G1},\cite{C}. The key idea is to replace the tangent bundle $T$ by $T\oplus T^*$ with its natural indefinite inner product $(X+\xi,X+\xi)=\iota_X\xi$ and the Lie bracket by the Courant bracket  
$$[X+\xi,Y+\eta]=[X,Y]+{\mathcal L}_X\eta-{\mathcal L}_Y\xi-\frac{1}{2}d(\iota_X\eta-\iota_Y\xi).$$
If $B$ is a closed $2$-form the action $X+\xi\mapsto X+\xi+\iota_XB$ preserves both the inner product and the Courant bracket and is called a B-field transform. 

A generalized complex structure is an orthogonal transformation $J:T\oplus T^*\rightarrow T\oplus T^*$ with $J^2=-1$ which satisfies an integrability condition which can be expressed in various ways, all analogous to the integrability condition for a complex structure but using the Courant bracket instead of the Lie bracket. The simplest is to take the isotropic subbundle $E$ of the complexification $(T\oplus T^*)^c$ on which $J=i$ and say that  sections of $E$ are closed under the Courant bracket. The standard examples are complex structures where $E$ is spanned by $(0,1)$ vector fields and $(1,0)$-forms, or symplectic structures where $E$ consists of objects of the form $X-i \iota_X\omega$ where $X$ is a vector field and $\omega$ the symplectic form. A holomorphic Poisson manifold defines a generalized complex structure where $E$ is spanned by $(0,1)$ vector fields  and objects of the form $ \sigma(\alpha)+\alpha$ where $\alpha$ is a $(1,0)$-form.

One of the key aspects of generalized geometry is that differential forms are interpreted as spinors -- the Clifford action of $T\oplus T^*$ on the exterior algebra of forms $\Lambda^*$ is 
$(X+\xi)\cdot \varphi=\iota_X\varphi+\xi\wedge \varphi$. Then the action of a 2-form $B$ is $\varphi\mapsto e^{-B}\varphi$ using exterior multiplication. There is an invariant pairing, the Mukai pairing, on forms with values in the top degree defined by $\langle \varphi,\psi \rangle =[\varphi\wedge s(\psi)]_n$ where $s(\psi)=\psi_0-\psi_1+\psi_2-\dots$, expanding by degree. 

Generalized complex structures are defined by maximal isotropic subbundles $E\subset (T\oplus T^*)^c$ and the annihilator under Clifford multiplication of any spinor is isotropic. If a complex form $\rho$ is closed and its annihilator is maximal isotropic (i.e. it is a pure spinor) with $E\cap \bar E=0$ (equivalently $\langle \rho,\bar\rho\rangle\ne 0$) then $\rho$ defines a generalized complex structure. An example is a symplectic structure where $\rho=e^{i\omega}$. The more general condition is that $E$ is integrable if 
$$d\rho=(X+\xi)\cdot\rho$$
for some local section $X+\xi$ of $(T\oplus T^*)^c$.

\subsection{Generalized Dolbeault operators}
If $f$ is a function on a generalized complex manifold $(M,J)$ we have 
$$df\in T^*\in (T\oplus T^*)^c=E\oplus \bar E$$
and we define $\bar\partial_J f$ to be the $\bar E$ component. For a complex structure this is the usual $\bar\partial f$ and for  a symplectic structure $\omega$, $\bar \partial_{J} f=(iX+df)/2$ where $X$ is the Hamiltonian vector field of $f$.  For a holomorphic Poisson structure $\sigma$ we obtain (where in the formula we use the standard meaning of $\bar\partial f$ and  $\partial f$):
\begin{equation}
\bar\partial_{J}f=\bar\partial f-\sigma(\partial f)+\bar\sigma(\bar\partial f)
\label{bard}
\end{equation}

The $\bar\partial_J$ operator can be extended to a generalized Dolbeault complex 
$$\cdots \rightarrow C^{\infty}(\Lambda^p\bar E)\stackrel{\bar\partial} \rightarrow C^{\infty}(\Lambda^{p+1}\bar E)\rightarrow\cdots$$
(where for simplicity we suppress the subscript). This is  purely analogous to the usual Dolbeault operator and it is well-defined because sections of the  bundle $E$ are closed under Courant bracket and $E$ is isotropic. It forms a complex for the same reason:  the term 
$([A,B],C) + ([B,C],A) + ([C,A],B)$,  
whose derivative obstructs the Jacobi identity, vanishes.

This motivates the definition of a {\it generalized holomorphic structure} on a vector bundle $V$ over a generalized complex manifold. This consists of a differential operator 
$$\bar\partial_V:C^{\infty}(V)\rightarrow C^{\infty}(V\otimes \bar E)$$
with the properties
\begin{itemize}
\item
$\bar\partial_V(fs)=\bar\partial f s+f \bar\partial_V s$
\item
$\bar\partial_V^2=0$
\end{itemize}
where the last condition involves the  bundle-valued extension of the generalized Dolbeault operator. If in a local basis the operator is defined by a matrix valued section $A$ of $\bar E$, then this condition is 
$\bar\partial A+A\cdot A=0$. (Note that since $\bar E$ is isotropic, for $e,e'\in \bar E$, $e\cdot e'=-e'\cdot e$ so that  this is essentially an exterior product $\bar\partial A+A\wedge A=0$.)

The Dolbeault complex is also related to the decomposition of forms on a generalized complex manifold.  The endomorphism $J$ of $T\oplus T^*$ is skew adjoint and we can consider its Lie algebra action on spinors (which of course are differential forms). If the manifold has (real) dimension $2m$, then the forms are decomposed into eigenbundles with eigenvalues $ik$ $(-m\le k\le m)$.
$$U_{-m},U_{-m+1}\dots,U_0,U_1,\dots,U_m.$$
\begin{ex}
 For a complex structure 
$$U_k=\bigoplus_{p-q=k}\Lambda^{p,q}.$$
\end{ex}

The integrability of the generalized complex structure means that the exterior derivative $d$ maps sections of $U_k$ to $U_{k-1}\oplus U_{k+1}$. The two parts are closely related to the $\bar\partial$  operators above, but we need to consider in more detail one of these eigenbundles first, namely  $U_m$.
\subsection{The canonical bundle}

The maximal isotropic subbundle $E\subset (T\oplus T^*)^c$ is the annihilator of a spinor, but a spinor only defined up to a scalar multiple so this defines a distinguished line bundle in $\Lambda^*$ called the {\it canonical bundle} $K$.

\begin{ex} On a complex $m$-dimensional manifold, $dz_1\wedge dz_2\wedge\dots \wedge dz_m$ is annihilated by interior product with a $(0,1)$ vector and exterior product with a $(1,0)$ form, so $K$ is the usual canonical bundle of a complex manifold. For a symplectic manifold, $e^{i\omega}$ trivializes the canonical bundle.
\end{ex}

In the eigenspace decomposition, $K$ is the subbundle of forms $U_m$. Moreover,
$$U_{m-k}\cong K\otimes \Lambda^k\bar E,$$
essentially generated by  the Clifford products of $k$ elements of $\bar E$ acting on $K$. 
Now, as remarked above, $d$ maps sections of $U_k$ to $U_{k-1}\oplus U_{k+1}$, and so takes sections of $U_m=K$ to $U_{m-1}=K\otimes \bar E$ since $U_{m+1}=0$. This defines a generalized holomorphic structure   on $K$. 

\begin{rmk} For a symplectic manifold, $e^{i\omega}$ trivializes $K$ and is closed and hence holomorphic in this generalized sense -- hence the appropriate  terminology {\it generalized Calabi-Yau manifold} for such a manifold.
\end{rmk}

When the canonical bundle is an even form there is a tautological section $\tau$ of its dual bundle $K^*$. This is just the projection from $\Lambda^*$ to $\Lambda^0=\C$, restricted to $K$. It is holomorphic in the generalized sense.  The section $\tau$ may be identically zero, but it is non-zero clearly at points where the generalized complex structure is the B-field transform of a symplectic structure, for 
$$e^Be^{i\omega}=1+(B+i\omega)+\dots$$ 

\begin{ex} For a holomorphic Poisson structure $\sigma$ on a surface, where $\sigma\ne 0$ the generalized complex structure is the B-field transform of a symplectic structure  ($\sigma^{-1} =B+i \omega$). The tautological  section vanishes on the elliptic curve $C$.
\end{ex}

Returning to  $d:C^{\infty}(U_k)\rightarrow C^{\infty}(U_{k-1}\oplus U_{k+1})$ we  write the projection to $U_{k-1}$ as $\bar\partial$ and to $U_{k+1}$ as $\partial$. The notation is consistent with the previous one in the sense that  $U_{m-1}=K\otimes \bar E$ and the operator is the $\bar\partial_K$ operator for the tautological generalized holomorphic structure on $K$. 

We have then a natural elliptic complex which we can write as either
$$\cdots \rightarrow C^{\infty}(K\otimes\Lambda^p\bar E)\stackrel{\bar\partial} \rightarrow C^{\infty}(K\otimes \Lambda^{p+1}\bar E)\rightarrow\cdots$$
or
$$\cdots \rightarrow C^{\infty}(U_{m-p})\stackrel{\bar\partial} \rightarrow C^{\infty}(U_{m-p-1})\rightarrow\cdots$$
\section{A generalized construction }
\subsection{The problem}
Suppose now that $M$ is a 4-manifold with a generalized complex structure such that the tautological section $\tau$ of the canonical bundle has a connected nondegenerate zero-set. As shown in \cite{CG1} this is a 2-torus with a complex structure, hence an elliptic curve $C$. We shall construct a rank 2 bundle on $M$ with a generalized holomorphic structure, given a set of points on $C$. 

We imitate the analytical approach to the Serre construction and take the bundle $1\oplus K^*$ where $K$ is the canonical bundle and find a distributional section $A$ of $K\otimes \bar E$ to define a generalized holomorphic structure by
$$\bar\partial (u,v)=(\bar\partial u+Av,\bar\partial v).$$
We are in the case $m=2$, so  $K\otimes \bar E\cong U_{2-1}=U_1$. We start with a set of points $x_i$ and look at the distributional form 
$$T=\sum_i\lambda_i\delta(x_i).$$
If we are to solve $\bar\partial A=T$ then $T$  must take values in $U_0$. There are then two questions that need to be answered:
\begin{enumerate}
\item
When does $T$ lie in $U_0$?
\item
When is $T=\bar\partial A$ for $A$ in $U_1$?
\end{enumerate} 
The first question needs a little more generalized geometry.

\subsection{Generalized complex submanifolds}

Given a submanifold $Y\subset M$, there is a distinguished subbundle 
$$TY\oplus N^*\subset (T\oplus T^*)\vert_Y$$
where $N^*$ is the conormal bundle. A submanifold is called {\it generalized complex} if $TY\oplus N^*$ is preserved by $J$. 
\begin{ex} For a complex manifold, this gives the usual notion of complex submanifold, for a symplectic manifold a Lagrangian submanifold. Applying a B-field to a symplectic structure can give new types of generalized complex submanifold but a point is never a generalized complex submanifold. Indeed a point $x$ is complex if the cotangent space $T^*_x\subset (T\oplus T^*)_x$ is preserved by $J$. But that means  there are complex cotangent vectors in $E$. However, $E$ is spanned by terms $X+\iota_X(B+i\omega)$ and so $X$ is never zero. 
\end{ex}
Now a compact oriented submanifold $Y^k$ defines a current $\Delta_Y$  in $\Omega^{n-k}$ by
$$\langle \Delta_Y, \alpha\rangle =\int_Y\alpha.$$
We then have
\begin{prp} $\Delta_Y$ lies in $U_0$ if and only if $Y$ is a generalized complex submanifold.
\end{prp}
\begin{prf} Consider the top exterior power $\Lambda^{2m-k}N^*$. Since $N^*\subset T$ is the annihilator of $TY$, if $\nu\in \Lambda^{2m-k}N^*$ is a generator, then $\iota_X\nu=0$ if and only if $X\in TY$. Similarly $\xi\wedge \nu=0$ if and only if $\xi\in N^*$. Thus $TY\oplus N^*$ is the annihilator under Clifford multiplication of $\nu$. 

If $Y$ is a generalized complex submanifold, then this annihilator is $J$-invariant, which means that  the real form $\nu$ is in the zero eigenspace of the Lie algebra action of $J$, i.e. $\nu \in U_0$. Conversely if $\nu\in U_0$, $Y$ is complex. 

Now consider a form $\alpha$.  To evaluate $\Delta_Y$ on this we take the degree $k$ component and integrate over $Y$. Now $\Lambda^{2m}T^*$ is canonically $\Lambda^{k}T^*Y\otimes \Lambda^{2m-k}N^*$. The Mukai pairing takes values in $\Lambda^{2m}T^*$, so  $\nu\mapsto\langle \alpha,\nu\rangle$ defines a homomorphism from $\Lambda^{2m-k}N^*$ to $\Lambda^{k}T^*Y\otimes \Lambda^{2m-k}N^*$, or equivalently an element of $\Lambda^{k}T^*Y$. It is straightforward to see that, up to a sign, this is the degree $k$ component of $\alpha$ restricted to $Y$. 

If $J\alpha=ik \alpha$, then 
$$ik\langle \alpha,\nu\rangle=\langle J\alpha,\nu\rangle=-\langle \alpha,J\nu\rangle=0$$

   Hence $\Delta_Y$ evaluated on $U_k$ for $k\ne 0$ is zero, and hence $\Delta_Y$ lies in $U_0$.
   \end{prf}   
   
   \begin{ex} The current defined by a complex submanifold of a complex manifold is of type $(p,p)$.
   \end{ex}
   
   Returning to our question we see that $T$ lies in $U_0$ if and only if each point $x_i$ is a generalized complex submanifold. Outside the elliptic curve $C$ the generalized complex structure is the B-field transform of a symplectic one, and as we have seen, points here are not complex.  In four dimensions, if $\tau=0$, the generalized complex structure is the stabilizer of a spinor of the form
   $e^B\alpha_1\wedge\alpha_2$. This is the B-field transform of an ordinary complex structure. Since $T^*$ is preserved by $J$ for an ordinary complex structure, and the B-field acts trivially on $T^*$ we see that any point on $C$ is a generalized complex submanifold. So we have an answer to the first question:
   \begin{prp} $T$ lies in $U_0$ if and only if the points $x_i$ lie on the elliptic curve $C$.
   \end{prp}
\subsection{The construction}
We now address the second question: suppose $T$ lies in $U_0$, when is it of the form $\bar\partial A$ for $A$ in $U_1$? In the standard case we used Serre duality to say that the Dolbeault cohomology class of $T$ is trivial if we evaluate on the generator $1$ of $H^0(M,{\mathcal O})$. 
 For the generalized $\bar\partial$ operator Serre duality consists of the non-degeneracy of the natural Mukai  pairing of $U_k$ and $U_{-k}$ at the level of cohomology.  A proof
 can be found in \cite{C}. 
 
 In our case, it means that $T$ in $U_0$ is cohomologically trivial if evaluation on all $\bar\partial$-closed forms in $U_0$ is zero. So suppose $\alpha$ is a section of $U_0$ with $\bar\partial\alpha=0$. The distribution $T$ is a sum of delta functions of points $x_i$, which lie on the curve $C$, so we need to know $U_0$ here.  But the generalized complex structure on $C$ is, as we have seen, the B-field transform of a complex structure. Now for a complex structure,
 $$U_0=\bigoplus_p \Lambda^{p,p}$$
 so
 $U_0\vert_C=e^B(\Lambda^{0,0}\oplus \Lambda^{1,1}\oplus \Lambda^{2,2})$, where $B$ is possibly locally defined.  Hence we can locally write 
 $$\alpha=e^B(a_0,a_1,a_2).$$
 However, $B$ leaves the degree zero part invariant, so
 in this local expression
 $a_0$ is the restriction of a globally defined function on $C$. 

Now, as shown in \cite{CG2}, a normal form (up to diffeomorphism and B-field transform) for a neighbourhood of a nondegenerate complex locus in four dimensions is provided by the holomorphic Poisson structure 
$$\sigma=z_1\frac{\partial}{\partial z_1}\wedge \frac{\partial}{\partial z_2}.$$
 
From this and (\ref{bard}) one can see that $\bar\partial \alpha=0$ implies that the degree zero term $a_0$ is holomorphic on the compact elliptic curve $C$, and hence constant.

Now $T$ evaluates at points $x_i\in C$ and involves just the degree zero component of $\alpha$. It follows that 
the condition on $T$ to be cohomologically trivial is 
$$\langle T,\alpha\rangle= const. \sum_i\lambda_i=0$$
as before.

We conclude:

\begin{thm} \label{main}  Let $M$ be a 4-manifold with a generalized complex structure such that the tautological section $\tau$ of the canonical bundle has a connected nondegenerate zero-set $C$. A set of $k$ distinct points $x_i\in C$ and $\lambda_i\in \C^*$ with $\lambda_1+\dots+\lambda_k=0$ defines a rank 2 generalized holomorphic bundle $V$ with a generalized holomorphic section  vanishing at the points $x_i$.
\end{thm}

\begin{rmk} Note that here we have no condition for a second section, but neither have we attempted to find bundles with {\it two} generalized holomorphic sections: the construction in Section \ref{con1} was a simple way to get Poisson modules, but they are more special than they need to be. So Theorem \ref{main} tells us that the Serre construction where the points are taken on a smooth anticanonical divisor gives us a Poisson module -- we don't need two sections of the line bundle on the elliptic curve. What happens is that the flat connection on $M\setminus C$ has upper-triangular  rather than trivial holonomy.\end{rmk}

\section{An application}
It is observed in  \cite{P} that if $V$ is a rank 2 Poisson module on a complex Poisson manifold, then the projective bundle $\PP(V)$ acquires a naturally induced Poisson structure. Here we prove a generalized version:
\begin{prp} Let $V$ be a rank two generalized holomorphic bundle over a generalized complex manifold $M$. Then $\PP(V)$ has a natural generalized complex structure.
\end{prp}
\begin{prf} First consider the generalized complex structure which is the product of the standard  complex structure on $\PP^1$ and the given generalized complex structure on $M$. If $\rho$ is a local non-zero section of the canonical bundle of $M$ then, using an affine coordinate $z$ on  $\PP^1$, $dz\wedge \rho$ is a section of the canonical bundle for the structure on the product. 

Now over an open set $U\subseteq M$ apply the  diffeomorphism of $\PP^1\times U$ defined by a map $a$ from $U$ to $SL(2,\C)$:
$$\tilde z=\frac{a_{11}z+a_{12}}{a_{21}z+a_{22}}.$$
Then
$$(a_{21}z+a_{22})^2d\tilde z= dz+A_{12} +(A_{11}-A_{22})z-A_{21}z^2=dz+\theta$$
where $A=a^{-1}da$.  Using $dA+A^2=0$, we have 
\begin{equation}
d(dz+\theta)=(dz+\theta)\wedge(2z A_{21}-(A_{11}-A_{22}))=(dz+\theta)\wedge\alpha.
\label{dtheta}
\end{equation}
Consider $(dz+\theta)\wedge \rho$, or equivalently the Clifford product $(dz+\theta)\cdot \rho$, since $dz+\theta$ is a one-form. Now $\rho$ is annihilated  by $E\subset (T\oplus T^*)^c$, so it is only the $\bar E$ component of $\theta$ (denote it $\theta^{01}$) which contributes. This defines an integrable generalized complex structure trivially since we simply transformed the product by a diffeomorphism, but a direct check of integrability goes as follows: we need to show that locally 
$$d((dz+\theta^{01})\cdot \rho)=\beta\cdot(dz+\theta^{01})\cdot\rho$$
for some section $\beta$ of $(T\oplus T^*)^c$.
But from (\ref{dtheta})
$$d((dz+\theta^{01})\wedge \rho)=-\alpha\wedge(dz+ \theta^{01})\wedge \rho-\theta^{01}\wedge\gamma \cdot\rho=-\alpha\cdot (dz+\theta^{01})\cdot \rho-(dz+\theta^{01})\cdot\gamma \cdot\rho$$
using the integrability $d\rho=\gamma\cdot\rho$ of the structure on $M$, where again we can take $\gamma$ to be in $\bar E$. Now since $\bar E$ is isotropic,  two sections anticommute under  the Clifford product , so
$$d((dz+\theta^{01})\wedge \rho)=(\gamma-\alpha)\cdot (dz+\theta^{01})\wedge \rho$$
which is the required integrability.

Now suppose $V$ is a rank $2$ bundle with a generalized holomorphic structure, and in a local trivialization $\bar\partial_V$ is defined by a ``connection matrix" $A$ with values in $\bar E$. Then define an almost  generalized  complex structure by 
$$(dz+A_{12} +(A_{11}-A_{22})z-A_{21}z^2)\cdot \rho.$$
In the argument for integrability above, we only needed the vanishing of the $\Lambda^2\bar E$ component of $dA+A^2=0$ and from the definition of a generalized holomorphic structure, we have  $\bar\partial A+A\cdot A=0$, so this provides the ingredient to prove integrability for the generalized holomorphic structure. 
\end{prf}

As a consequence, if we use our construction to generate rank 2 vector bundles with generalized holomorphic structure on the Cavalcanti-Gualtieri 4-manifolds, we can find six-dimensional generalized complex examples on their projective bundles. These have a structure which is {\it complex} in the fibre directions. In \cite{C0} it is shown that a symplectic bundle over a generalized complex base has a  generalized complex structure which is {\it symplectic} along the fibres.

\vskip 1cm
 Mathematical Institute, 24-29 St Giles, Oxford OX1 3LB, UK
 
 hitchin@maths.ox.ac.uk

 \end{document}